\documentclass[reqno,12pt]{amsart}

\usepackage{amsmath,amsfonts,amsthm,amscd,upref,amstext}
\usepackage{amssymb}
\usepackage{mathrsfs}
\usepackage{mathtools}
\usepackage{xcolor}
\usepackage{fullpage}
\usepackage{lmodern}

\theoremstyle{plain}
\newtheorem{prop}{Proposition}[section]
\newtheorem{thm}[prop]{Theorem}
\newtheorem{cor}[prop]{Corollary}
\newtheorem{lem}[prop]{Lemma}

\theoremstyle{definition}
\newtheorem{defn}[prop]{Definition}
\newtheorem{ex}[prop]{Example}
\newtheorem{question}[prop]{Question}
\newtheorem{nota}[prop]{Notation}

\theoremstyle{remark}
\newtheorem{rem}[prop]{Remark}

\newcommand{\m}{\mathfrak m}
\newcommand{\cR}{\mathcal R}
\newcommand{\cG}{\mathcal G}
\DeclareMathOperator{\Spec}{Spec}
\DeclareMathOperator{\Proj}{Proj}

\DeclareMathOperator{\gr}{gr}
\DeclareMathOperator{\rank}{rank}
\DeclareMathOperator{\Sym}{Sym}
\DeclareMathOperator{\br}{br}
\DeclareMathOperator{\Hilb}{Hilb}

\numberwithin{equation}{section}

\begin{document}

\title{On a Question of Cid-Ruiz: positive results and counterexamples}
\date{}

\author{T. H. Freitas}
\address{Universidade Tecnol\'ogica Federal do Paran\'a, 85053-525,
Guarapuava-PR, Brazil}
\email{freitas.thf@gmail.com}

\author{V. H. Jorge P\'erez}
\address{Universidade de S{\~a}o Paulo -- ICMC, 13560-970,
S{\~a}o Carlos-SP, Brazil}
\email{vhjperez@icmc.usp.br}


\begin{abstract}
In this paper we give positive criteria and counterexamples to a question of
Cid-Ruiz on whether intrinsic polar multiplicities detect reductions
of modules. Over an excellent regular local surface with infinite
residue field, we compute the entire polar sequence of every
finite-colength submodule of a free module, showing that its only
variable entry is the order of the determinant ideal. This computation
yields classes for which equality of polar multiplicities
characterizes reduction, while also producing counterexamples in
every positive rank, including integrally closed modules with
identical intrinsic Hilbert functions and arbitrarily large relative
Buchsbaum--Rim multiplicity. In arbitrary dimension, we prove that,
under equidimensionality and a sharp dimension bound on the kernel of
the induced map of associated graded Rees algebras, equality of the
sums of the polar multiplicities characterizes reduction.
\end{abstract}

\maketitle

\section{Introduction}
\label{sec:introduction}

Multiplicity theory provides a natural way to study integral
dependence through numerical invariants. A classical example is
Rees's criterion for primary ideals, whose extension to modules
draws on Buchsbaum--Rim multiplicity and the theory of reductions
\cite{BuchsbaumRim1964,Rees1987,KT1994,SUV2001}. More recently,
mixed multiplicities, epsilon multiplicity, and Segre classes have
led to further criteria for integral dependence
\cite{DRTideals,DRTmodules,CidRuizSegre}. In this direction,
Cid-Ruiz \cite{CidRuiz2024} uses the polar multiplicities of
Kleiman--Thorup to study integral dependence of graded algebras and
modules.

Let $(R,\m)$ be an excellent local domain of dimension $d\geq1$.
For a finite torsion-free module $M$ of rank $e\geq1$, write
$\cR(M)=\bigoplus_{v\geq0}M^v$ and $r=d+e-1$. The intrinsic polar
multiplicities $\br_0(M),\ldots,\br_r(M)$ are the normalized leading
coefficients of the eventual polynomial associated with
$H_M(v,n)=\lambda_R(M^v/\m^{n+1}M^v)$. Since these invariants are
preserved under reduction, Cid-Ruiz asked the following.

\begin{question}[Cid-Ruiz, {\cite[Question~5.2]{CidRuiz2024}}]
\label{question:main}
Under which conditions on $U\subseteq E$ does one have
$U$ is a reduction of $E$ if and only if
$\br_i(U)=\br_i(E)$ for all $i$?
\end{question}

The failure of polar multiplicities to detect integral dependence
for ideals already appears in \cite[Example~6.9]{CPU}. Our aim is to
describe more precisely both sides of the problem, i.e., situations in
which the converse in Question~\ref{question:main} holds, and how
strongly it may fail otherwise.

Our first main result gives a complete description on regular
surfaces. Let $(R,\m,k)$ be an excellent regular local ring of
dimension two with infinite residue field, let $F=R^e$, and let
$M\subseteq F$ have finite colength. If $\mathfrak d_F(M)$ denotes
the determinant ideal of $M$, then
Theorem~\ref{thm:surface-polars} gives that
$\br_1(M)=\operatorname{ord}_{\m}\mathfrak d_F(M)$,
$\br_2(M)=1$, and $\br_i(M)=0$ for $i\neq1,2$.
In addition, for ideals, this result recovers the classical mixed-multiplicity formula
on a regular surface \cite{Verma1988,Trung2001}, while the theorem
applies to arbitrary finite-colength modules.

Together with Rees's determinantal criterion
\cite{Rees1987,HunekeSwanson}, this computation explains the main
obstruction in Question~\ref{question:main}: equality of polar
sequences records only the $\m$-adic order of the determinant ideal,
whereas reduction is governed by its integral closure. For instance, if
$\overline{\mathfrak d_F(U)}=\overline{Q^a}$ and
$\overline{\mathfrak d_F(E)}=\overline{Q^b}$ for powers of a fixed
$\m$-primary ideal $Q$, equality of the polar sequences
characterizes reduction
(Corollary~\ref{cor:fixed-determinant}). If
$\overline{\mathfrak d_F(U)}=\m^a$, equality of $\br_1$ alone
suffices (Corollary~\ref{cor:maximal-determinant}).

A second positive answer  holds in arbitrary dimension. For
$U\subseteq E$ of the same rank, consider
$\varphi:\gr_{\m}\cR(U)\longrightarrow\gr_{\m}\cR(E)$.
If $\gr_{\m}\cR(E)$ is equidimensional and
$\dim\ker\varphi<d+e$, then Theorem~\ref{thm:small-kernel} shows that
$U$ is a reduction of $E$ if and only if
$\sum_{i=0}^{r}\br_i(U)=\sum_{i=0}^{r}\br_i(E)$.
Thus, under these hypotheses, equality of the sums of the polar
multiplicities already detects reduction. The criterion applies, in
particular, to homogeneous modules generated in the same degree
(Corollary~\ref{cor:equigenerated}).

The surface formula also produces a large families of
counterexamples. Actually, Proposition~\ref{prop:closed-counterexamples} shows
that even when the smaller module is integrally closed, equality of
the polar sequences need not imply reduction. A more explicit
example is given by the ideals $I_q=(x,y^q)\subseteq k[[x,y]]$, where 
Proposition~\ref{prop:plane-example} shows that all of them have the
same intrinsic Hilbert function, although $I_q$ is not a reduction
of $(x,y)$ for $q\geq2$. Also, more  counterexamples are provided by the normal monomial ideals
$K_{p,q}=(x^2,xy^p,y^{p+q})$. By
Theorem~\ref{thm:two-face-family}, the corresponding modules
$K_{p,q}R^e$ have intrinsic Hilbert functions independent of the
parameters $p$ and $q$, while their relative behavior varies. In this sense, 
Corollary~\ref{cor:unbounded-two-face} gives, in every fixed positive
rank, proper inclusions $U\subsetneq E$ which are not reductions but
satisfy $H_U=H_E$, while $\br_0(U,E)$ can be arbitrarily large.
Thus even the entire intrinsic Hilbert functions of the two modules
do not control the relative defect of the inclusion.

It is important to realize that the dimension condition in Theorem~\ref{thm:small-kernel} is also
sharp. For the family $(x,y^q)\subsetneq(x,y)$,
Proposition~\ref{prop:explicit-kernel} shows that the kernel of the
associated graded map has dimension exactly $d+e$. Furthermore, we study the relative invariants that retain information
about the inclusion itself. For finite-length pairs,
$\br_0(U,E)$ agrees with the relative Buchsbaum--Rim multiplicity,
and Theorem~\ref{thm:twist-polynomial} expresses
$e_{\mathrm{pair}}(\m^aM,M)$ directly in terms of the absolute polar
multiplicities of $M$.

The paper is organized as follows.
Section~\ref{sec:preliminaries} recalls Rees algebras, reductions,
and the basic properties of intrinsic polar multiplicities.
Section~\ref{sec:surfaces} treats regular surfaces and the
determinantal criteria, while Section~\ref{sec:criterion} proves the
associated graded criterion in arbitrary dimension.
Section~\ref{sec:relative} develops relative polar and
Buchsbaum--Rim multiplicities, and
Section~\ref{sec:ideal-multiples} studies ideal multiples and
maximal-ideal twists. The counterexamples, including the normal
monomial families and the sharpness of the kernel bound, are
given in Section~\ref{sec:counterexamples}.
\section{Rees algebras and preliminary properties of polar multiplicities}
\label{sec:preliminaries}

We begin by fixing the notation and recalling the basic facts on Rees
algebras, reductions, and polar multiplicities that will be used
throughout the paper. We also record some elementary consequences that
will play an important role later.
Throughout the paper, $(R,\m,k)$ is an excellent local domain
of dimension $d\geq1$, unless explicitly stated otherwise. Write
$Q=\operatorname{Frac}(R)$. For a finite torsion-free module $M$ of
rank $e\geq1$, its Rees algebra is
\[
\cR(M)=\Sym_R(M)/(R\text{-torsion})
       =\bigoplus_{v\geq0}M^v.
\]

We view it inside $\Sym_Q(M\otimes_R Q)$. In particular, it is a
domain of dimension $d+e$. For $U\subseteq E$ of the same rank, put $A=\cR(U)$,  $B=\cR(E)$, and  $r=d+e-1$.
These are standard graded $R$-algebras, with a natural homogeneous
inclusion $A\subseteq B$.

\begin{defn}\label{def:reduction}
The module $U$ is a \emph{reduction} of $E$ if
$B_{v+1}=A_1B_v$ for some $v\geq0$, or equivalently for all
sufficiently large $v$. Equivalently, $B$ is finite over $A$.
\end{defn}

For $M\subseteq F$ with $F$ free, its integral closure
$\overline M$ in $F$ consists of the elements whose degree-one
images in $\Sym_R(F)$ are integral over $\cR(M)$.
A module is \emph{integrally closed}, or \emph{complete}, if
$M=\overline M$. For $U\subseteq E\subseteq F$, reduction
is equivalent to $E\subseteq\overline U$ (see
\cite{Rees1987,HunekeSwanson}).

\begin{nota}\label{not:hilbert}
For $C=\cR(M)$, write $\cG_M=\gr_{\m C}(C)$,  and
$(\cG_M)_{(v,n)}=\m^nM^v/\m^{n+1}M^v.$
The function $H_M(v,n)=\lambda_R(M^v/\m^{n+1}M^v)$
agrees with a polynomial for $v,n\gg0$. By
\cite[Lemmas~2.3 and~2.4]{CidRuiz2024}, its top-degree part is given by 
\begin{equation}\label{eq:polar-polynomial}
\sum_{i=0}^r\frac{\br_i(M)}{(r-i)!i!}v^{r-i}n^i.
\end{equation}
When an associated graded ring is given a single grading below, we
use total degree $v+n$.
\end{nota}

\begin{defn}\label{def:geometric-polar}
Let $X_M=\Proj(\cR(M))$, let $p_M:X_M\to\Spec(R)$ be the
structure morphism, and let $Z_M=p_M^{-1}(\{\m\})$ with its
scheme structure. Write
\[
\mathscr D_M=\mathbb P(C_{Z_M}X_M\oplus\mathbf1),\qquad
\pi_M:\mathscr D_M\longrightarrow X_M,
\]
for the projective completion of the normal cone and its natural
projection. Then the geometric definition is
\begin{equation}\label{eq:geometric-polar}
\br_i(M)=\int_{\mathscr D_M}
c_1(\mathcal O_{\mathscr D_M}(1))^i
c_1(\pi_M^*\mathcal O_{X_M}(1))^{r-i}
\cap[\mathscr D_M]_r
=m_r^i(\cR(M)).
\end{equation}
Here $[\mathscr D_M]_r$ is the dimension-$r$ part of the
fundamental cycle. The equivalence with
\eqref{eq:polar-polynomial} is
\cite[Definition~2.2 and Lemma~2.4]{CidRuiz2024}.
The scheme $\mathscr D_M$ is distinct from the graded algebra
$D_M$ introduced in Lemma~\ref{lem:total-multiplicity}.
\end{defn}

\begin{prop}
\label{prop:fixed-entries}
Under the standing hypotheses, assume that $M$ has rank $e\geq1$, and
write $e_{\m}(R)$ for the Hilbert--Samuel multiplicity of $R$. Then
\begin{equation}\label{eq:fixed-entries}
\br_0(M)=0,\qquad \br_d(M)=e_{\m}(R),\qquad
\br_i(M)=0\quad(d<i\leq r).
\end{equation}
In particular, when the base ring and the rank are fixed, only the
entries $\br_i(M)$ with $1\leq i\leq d-1$ can depend on $M$.
\end{prop}

\begin{proof}
Since $\cR(M)$ is $R$-torsion-free and $R$ is a domain, every nonzero
element of $\m$ is a nonzerodivisor on $\cR(M)$. Hence
$H^0_{\m}(\cR(M))=0$. By
\cite[Proposition~2.10(i)]{CidRuiz2024}, the zeroth polar multiplicity
is the corresponding $j$-multiplicity, and therefore
$\br_0(M)=0$. The vanishing $\br_i(M)=0$ for $i>d$ follows from
\cite[Lemma~4.2]{CidRuiz2024}. For each $v\geq0$, the coefficient of $n^d$ in the Hilbert polynomial
of $H_M(v,n)=\lambda_R(M^v/\m^{n+1}M^v)$
is $e_{\m}(M^v)/d!$. Since $R$ is a domain and $\rank_R(M^v)=\binom{v+e-1}{e-1},$
the associativity formula for Hilbert--Samuel multiplicity gives $e_{\m}(M^v)
=
\binom{v+e-1}{e-1}e_{\m}(R).$
Thus the coefficient of $v^{e-1}n^d$ in the top-degree part of
$H_M(v,n)$ is $\frac{e_{\m}(R)}{(e-1)!d!}.$
On the other hand, since $r=d+e-1$, the term with $i=d$ in
\eqref{eq:polar-polynomial} contributes
\[
\frac{\br_d(M)}{(r-d)!d!}v^{r-d}n^d
=
\frac{\br_d(M)}{(e-1)!d!}v^{e-1}n^d.
\]
Comparing coefficients yields
$\br_d(M)=e_{\m}(R)$.
\end{proof}

\begin{cor}\label{cor:dimension-one}
Let $(R,\m)$ be a one-dimensional excellent local domain. Then all finite
torsion-free $R$-modules of a fixed positive rank $e$ have the same
polar sequence:
\[
(\br_0(M),\ldots,\br_e(M))=(0,e_{\m}(R),0,\ldots,0).
\]
\end{cor}

\begin{proof}
Apply Proposition~\ref{prop:fixed-entries} with $d=1$.
\end{proof}

\begin{prop}\label{prop:necessary}
If $U$ is a reduction of $E$, then $\br_i(U)=\br_i(E)$
for every $i$.
\end{prop}

\begin{proof}
Set $A=\cR(U)$ and $B=\cR(E)$. Since $U$ is a reduction of $E$,
the inclusion $A\subseteq B$ is finite. Moreover,
$U\otimes_R Q=E\otimes_R Q$, and hence $A\otimes_R Q=B\otimes_R Q
\simeq \Sym_Q(E\otimes_R Q),$
which is a polynomial ring over $Q$. Thus $A\subseteq B$ is a finite
birational extension. The conclusion now follows from
\cite[Theorem~4.4(iii),(iv)]{CidRuiz2024}.
\end{proof}

\section{Polar multiplicities on a regular surface}
\label{sec:surfaces}

Throughout this section, $(R,\m,k)$ is an excellent regular local
ring of dimension two with infinite residue field, and $F=R^e$.
For a finite-colength submodule $M\subseteq F$, write
\[
\mathfrak d_F(M)=I_e(M)=\operatorname{Fitt}_0(F/M),\qquad
o(M)=\operatorname{ord}_{\m}\mathfrak d_F(M).
\]
Here $I_e(M)$ is the ideal of maximal minors of a matrix of
generators of $M$, and $\operatorname{ord}_{\m}L$
is the largest $a$ such that $L\subseteq\m^a$. We set
$\mathfrak d_F(F)=R$ and $o(F)=0$. The ambient module will be
omitted from the notation when it is fixed.

The following formula describes the entire intrinsic polar sequence.
It does not require integral closedness, a homogeneous presentation,
or a decomposition into ideals.

\begin{thm}\label{thm:surface-polars}
For every finite-colength submodule $M\subseteq F$,
\begin{equation}\label{eq:surface-polars}
\br_1(M)=o(M),\qquad \br_2(M)=1,\qquad
\br_i(M)=0\quad(i=0\text{ or }3\leq i\leq e+1).
\end{equation}
\end{thm}

To prove Theorem~\ref{thm:surface-polars}, we first establish two
auxiliary lemmas concerning Hilbert coefficients and determinant
orders of Rees components.

\begin{lem}\label{lem:surface-hilbert-coefficient}
Let $N\subseteq R^s$ have finite colength, and write
\[
\lambda_R(N/\m^{n+1}N)
=s\binom{n+2}{2}-e_1^{\m}(N)(n+1)+e_2^{\m}(N)
\quad(n\gg0).
\]
Then $e_1^{\m}(N)=-\operatorname{ord}_{\m}I_s(N).$
\end{lem}

\begin{proof}
Set $L=R^s/N$ and $b=\operatorname{ord}_{\m}I_s(N)$. Since the residue
field is infinite, we may choose $x\in\m\setminus\m^2$ superficial for
the $\m$-adic filtration of $N$ such that the image of $I_s(N)$ in
$S=R/(x)$ still has order $b$. Put $T=(0:_Lx)$, and  $\overline N=\operatorname{im}(N\to S^s).$
Tensoring $0\to N\to R^s\to L\to0$ with $S$ gives
\[
0\longrightarrow T\longrightarrow N/xN\longrightarrow\overline N
\longrightarrow0.
\]
Since $S$ is a DVR and $\overline N$ is torsion-free of rank $s$, it is
free. Hence, for $n\gg0$,
\[
\lambda_S\bigl((N/xN)/\m^{n+1}(N/xN)\bigr)
=s(n+1)+\lambda_R(T).
\]
On the other hand, superficiality of $x$ gives
\[
\lambda_S\bigl((N/xN)/\m^{n+1}(N/xN)\bigr)
=s(n+1)-e_1^{\m}(N)
\]
for $n\gg0$. Therefore $-e_1^{\m}(N)=\lambda_R(T)$.
Since $L$ has finite length, multiplication by $x$ on $L$ has kernel
and cokernel of the same length, so $\lambda_R(T)=\lambda_R(L/xL)
=\lambda_S(S^s/\overline N).$
Over the DVR $S$, the latter length is the valuation of the determinant
ideal $I_s(\overline N)$. By the choice of $x$, this valuation is $b$.
Thus $-e_1^{\m}(N)=b=\operatorname{ord}_{\m}I_s(N),$
which proves the result.
\end{proof}

\begin{lem}\label{lem:symmetric-determinant-order}
For every integer $v\geq1$,
\begin{equation}\label{eq:symmetric-determinant-order}
\operatorname{ord}_{\m}
\mathfrak d_{\Sym_R^v(F)}(M^v)
=
\binom{v+e-1}{e}\,o(M).
\end{equation}
\end{lem}

\begin{proof}
Since $\gr_{\m}R$ is a domain, the $\m$-adic order extends to a
discrete valuation $\nu$ of $Q=\operatorname{Frac}(R)$, given by $\nu(a/b)=\operatorname{ord}_{\m}(a)-\operatorname{ord}_{\m}(b).$
Let $V$ be its valuation ring and let $t$ be a uniformizer. By the
elementary divisor theorem, after choosing suitable bases of
$F\otimes_RV$, we may write
\[
MV=\bigoplus_{j=1}^e t^{a_j}Vf_j,
\qquad a_j\geq0.
\]
Since the determinant ideal commutes with extension to $V$, we have $\sum_{j=1}^e a_j=o(M).$

The extension of $M^v$ to $V$ has basis $t^{\sum_{j=1}^e a_j\alpha_j}
f_1^{\alpha_1}\cdots f_e^{\alpha_e},$ where
$|\alpha|=v.$
Hence the valuation of its determinant is
\[
\sum_{|\alpha|=v}\sum_{j=1}^e a_j\alpha_j
=
\sum_{j=1}^e a_j
\sum_{|\alpha|=v}\alpha_j.
\]
For each $j$,
\[
\sum_{|\alpha|=v}\alpha_j
=
\binom{v+e-1}{e},
\]
and therefore
\[
\nu\!\left(
\mathfrak d_{\Sym_V^v(F\otimes_RV)}(M^vV)
\right)
=
\binom{v+e-1}{e}\sum_{j=1}^e a_j
=
\binom{v+e-1}{e}\,o(M).
\]
Since maximal minors commute with extension to $V$, this is precisely
\eqref{eq:symmetric-determinant-order}.
\end{proof}
\begin{proof}[Proof of Theorem~\ref{thm:surface-polars}]
For every $v\geq1$, the module $M^v$ has finite colength in
$\Sym_R^v(F)$, since $M_{\mathfrak p}=F_{\mathfrak p}$ for every
$\mathfrak p\neq\m$. Applying Lemmas~\ref{lem:surface-hilbert-coefficient}
and~\ref{lem:symmetric-determinant-order} gives, for $n\gg0$,
\begin{equation}\label{eq:surface-hilbert-polynomial}
H_M(v,n)
=
\binom{v+e-1}{e-1}\binom{n+2}{2}
+
\binom{v+e-1}{e}o(M)(n+1)
+
e_2^{\m}(M^v).
\end{equation}

For $v\gg0$, this agrees with the bivariate Hilbert polynomial in
Notation~\ref{not:hilbert}. Hence the coefficient of $n$ is
\[
\frac32\binom{v+e-1}{e-1}
+
\binom{v+e-1}{e}o(M).
\]
Its degree-$e$ term in $v$ is $\frac{o(M)}{e!}v^e.$
Since $d=2$ and hence $r=e+1$, the term with $i=1$ in
\eqref{eq:polar-polynomial} is $\frac{\br_1(M)}{e!}v^e n.$
Comparing coefficients yields $\br_1(M)=o(M).$
Finally, Proposition~\ref{prop:fixed-entries} gives
$\br_2(M)=1$ and
$\br_i(M)=0$ for $i=0$ or $3\leq i\leq e+1$.
\end{proof}

\begin{rem}\label{rem:surface-literature}
For an $\m$-primary ideal $I$, the identity $H_I(v,n)=\lambda_R(R/\m^{n+1}I^v)-\lambda_R(R/I^v)$
identifies $\br_1(I)$ with the mixed multiplicity
$e(\m^{[1]},I^{[1]})$. On a regular surface this equals
$\operatorname{ord}_{\m}I$ (see
\cite[Corollary~4.3]{Trung2001} and
\cite[proof of Theorem~4.1]{Verma1988}).
Thus Theorem~\ref{thm:surface-polars} extends this formula from ideals
to arbitrary finite-colength modules through their determinant ideals.
\end{rem}


\begin{prop}
\label{prop:determinant-reduction}
Let $R$ be a Noetherian domain and let $U\subseteq E\subseteq R^e$
have rank $e$. Then $U$ is a reduction of $E$ if and only if
$\overline{\mathfrak d(U)}
=
\overline{\mathfrak d(E)}.$
\end{prop}

\begin{proof}
This is Rees's determinantal criterion (see
\cite{Rees1987} or \cite[Proposition~16.3.2]{HunekeSwanson}).
\end{proof}

\begin{cor}
\label{cor:fixed-determinant}
Let $Q$ be an $\m$-primary ideal of $R$, and let
$U\subseteq E\subseteq F$ have finite colength. Suppose that $\overline{\mathfrak d(U)}=\overline{Q^a}$, and  $\overline{\mathfrak d(E)}=\overline{Q^b}$
for some $a,b\geq1$. Then $U$ is a reduction of $E$ if and only if
their intrinsic polar sequences are equal.
\end{cor}

\begin{proof}
Since the $\m$-adic order is unchanged by integral closure, $o(U)=a\operatorname{ord}_{\m}(Q)$, and $o(E)=b\operatorname{ord}_{\m}(Q).$
As $\operatorname{ord}_{\m}(Q)>0$, equality of the polar sequences,
together with Theorem~\ref{thm:surface-polars}, implies $a=b$.
Hence $\overline{\mathfrak d(U)}
=
\overline{\mathfrak d(E)},$
and Proposition~\ref{prop:determinant-reduction} gives that $U$ is a
reduction of $E$. The converse follows from
Proposition~\ref{prop:necessary}.
\end{proof}

\begin{cor}\label{cor:maximal-determinant}
Let $U\subseteq E\subseteq F$ have finite colength and suppose that $\overline{\mathfrak d(U)}=\m^a$
for some $a\geq0$. Then $U$ is a reduction of $E$ if and only if $\br_1(U)=\br_1(E).$
\end{cor}

\begin{proof}
Suppose that $\br_1(U)=\br_1(E)$. By
Theorem~\ref{thm:surface-polars}, $o(E)=o(U)=a.$
Since $U\subseteq E$, we have $\mathfrak d(U)\subseteq\mathfrak d(E)\subseteq\m^a.$ As $\m^a$ is integrally closed, taking integral closures yields
\[
\m^a=\overline{\mathfrak d(U)}
\subseteq\overline{\mathfrak d(E)}
\subseteq\m^a.
\]
Hence
$\overline{\mathfrak d(U)}=\overline{\mathfrak d(E)}$, and
Proposition~\ref{prop:determinant-reduction} shows that $U$ is a
reduction of $E$. The converse follows from
Proposition~\ref{prop:necessary}.
\end{proof}

The preceding criteria involve only determinant ideals and impose no
homogeneity condition on the generators. The following example
illustrates this point.

\begin{ex}\label{ex:surface-positive}
Let $R=k[[x,y]]$, with $k$ infinite, and set $E=\m f_1\oplus\m^2 f_2$, and $U=R(x,0)+R(y,x^2)+R(0,y^2)\subseteq E.$
The maximal minors give $\mathfrak d(U)=(x^3,xy^2,y^3),$ and 
$\mathfrak d(E)=\m^3$. Since $(x^3,y^3)\subseteq\mathfrak d(U)\subseteq\m^3$
and $\overline{(x^3,y^3)}=\m^3$, it follows that $\overline{\mathfrak d(U)}=\m^3
=\overline{\mathfrak d(E)}.$
Hence $U$ is a reduction of $E$ by
Proposition~\ref{prop:determinant-reduction}. Moreover,
Theorem~\ref{thm:surface-polars} gives the common polar sequence $(0,3,1,0).$
The inclusion is proper: $U$ is generated by three elements, whereas $\mu(E)=\mu(\m)+\mu(\m^2)=5.$
Hence, $U$ has finite colength in $R^2$ because
$\mathfrak d(U)$ is $\m$-primary.
\end{ex}

\section{A criterion via associated graded rings}
\label{sec:criterion}

In this section, $e(-)$ denotes
the usual multiplicity of a standard graded algebra over a field.

\begin{lem}\label{lem:total-multiplicity}
Let $M$ be as in Section~\ref{sec:preliminaries}, of rank $e$, and set
$r=d+e-1$. Let $D_M=\cG_M[z]$, where $\deg z=(0,1)$, and endow
$D_M$ with the total grading. Then $D_M$ is a standard graded
$k$-algebra of dimension $r+2$, and $e(D_M)=\sum_{i=0}^r\br_i(M).$

\end{lem}

\begin{proof}
Since the bigraded component $(D_M)_{(v,n)}$ has dimension
$H_M(v,n)$, we have
\begin{equation}\label{eq:total-hilbert}
\dim_k(D_M)_N
=
\sum_{v=0}^N H_M(v,N-v).
\end{equation}

Choose $c$ such that the Hilbert polynomial for $H_M(v,n)$ is valid
whenever $v,n\geq c$. The terms with $v<c$ contribute $O(N^d)$.
For each fixed $n<c$, the function $H_M(v,n)$ is the Hilbert function
of the standard graded algebra
$\cR(M)/\m^{n+1}\cR(M)$, whose dimension is at most $r+1$; hence these
terms contribute $O(N^r)$. Since $d\leq r$, all boundary terms are
$O(N^r)$.

On the remaining range, substitute \eqref{eq:polar-polynomial}. The
terms of total degree $<r$ contribute $O(N^r)$, while for each
$0\leq i\leq r$,
\[
\sum_{v=0}^N
\frac{v^{r-i}(N-v)^i}{(r-i)!i!}
=
\frac{N^{r+1}}{(r+1)!}+O(N^r).
\]
Therefore
\[
\dim_k(D_M)_N
=
\frac{\sum_{i=0}^r\br_i(M)}{(r+1)!}N^{r+1}
+O(N^r).
\]
Since the top-degree part of $H_M(v,n)$ is nonzero, the sum of the
polar multiplicities is positive. Hence $\dim D_M=r+2$, and comparison
with the Hilbert polynomial of the standard graded algebra $D_M$ gives $e(D_M)=\sum_{i=0}^r\br_i(M).$
\end{proof}


\begin{lem}\label{lem:graded-criterion}
Let $C\subseteq D$ be a homogeneous inclusion of standard graded
algebras over a field, with $\dim C=\dim D$ and $D$
equidimensional. If $e(C)=e(D),$ then $D$ is finite over $C$.
\end{lem}

\begin{proof}
This follows from \cite[Proposition~6.1(b)]{SUV2001} (see also the
discussion following Question~1.1 in \cite{CidRuiz2024}).
\end{proof}

\begin{nota}\label{not:kernel}
For $A=\cR(U)\subseteq B=\cR(E)$, let $\varphi:\cG_U\longrightarrow\cG_E$ and $K=\ker\varphi$. Then $K$ is a bihomogeneous ideal of $\cG_U$. We write $\dim K$ for the Krull dimension of its support as a $\cG_U$-module, with the convention $\dim 0=-\infty$. Explicitly,
\begin{equation}\label{eq:kernel-components}
K_{(v,n)}=
\frac{\m^nA_v\cap\m^{n+1}B_v}{\m^{n+1}A_v}.
\end{equation}
\end{nota}

We now obtain a sufficient criterion in arbitrary dimension. Under a
small-kernel hypothesis on the associated graded map, equality of the
sum of the polar multiplicities already detects reduction.

\begin{thm}\label{thm:small-kernel}
Let $(R,\m)$ be an excellent local domain of dimension $d\geq1$, and
let $U\subseteq E$ be finite torsion-free modules of the same rank
$e\geq1$. Suppose that
\begin{enumerate}
\item $\cG_E=\gr_{\m}\cR(E)$ is equidimensional;
\item $\dim K<d+e$, where $K$ is as in
Notation~\ref{not:kernel}.
\end{enumerate}
Then the following conditions are equivalent:
\begin{enumerate}
\item[(i)] $U$ is a reduction of $E$;
\item[(ii)] $\br_i(U)=\br_i(E)$ for all $0\leq i\leq r$;
\item[(iii)] $\sum_{i=0}^r\br_i(U)=\sum_{i=0}^r\br_i(E)$.
\end{enumerate}
\end{thm}

\begin{proof}
By Proposition~\ref{prop:necessary}, (i)$\Rightarrow$(ii), while
(ii)$\Rightarrow$(iii) is immediate. 

Suppose (iii). Extend $\varphi:\cG_U\to\cG_E$ by the identity on $z$ to a homogeneous
homomorphism $D_U\to D_E$, with notation as in
Lemma~\ref{lem:total-multiplicity}, and let $C$ be its image. Then
\[
0\longrightarrow K[z]\longrightarrow D_U
\longrightarrow C\longrightarrow0
\]
is an exact sequence of graded $D_U$-modules. Since
$\dim D_U=d+e+1$ and
\[
\dim K[z]=\dim K+1<d+e+1,
\]
the top Hilbert coefficient is unchanged, and hence $\dim C=d+e+1$, and  
$e(C)=e(D_U).$
By Lemma~\ref{lem:total-multiplicity} and (iii),  one has $e(C)=e(D_U)=e(D_E).$
Moreover, $D_E$ is equidimensional because $\cG_E$ is
equidimensional. Thus Lemma~\ref{lem:graded-criterion}, applied to
$C\subseteq D_E$, shows that $D_E$ is finite over $C$.

Passing first modulo $z$ and then modulo the ideal generated by the
positive-degree part of $\gr_{\m}R$, we obtain a finite homomorphism
of special fiber rings $A/\m A\longrightarrow B/\m B.$
Since these are standard graded in Rees degree, finiteness implies
that, for $v\gg0$,
\[
B_{v+1}/\m B_{v+1}
=
A_1\bigl(B_v/\m B_v\bigr).
\]
Equivalently, $B_{v+1}=A_1B_v+\m B_{v+1}.$
The quotient $B_{v+1}/A_1B_v$ is a finite $R$-module, so Nakayama's
lemma gives $B_{v+1}=A_1B_v$
 for $v\gg0$.
Hence $U$ is a reduction of $E$, proving (iii)$\Rightarrow$(i).
\end{proof}

\begin{lem}\label{lem:strictness}
The map $\varphi$ is injective if and only if
\begin{equation}\label{eq:strictness}
A_v\cap\m^nB_v=\m^nA_v
\qquad\text{for all }v,n\geq0.
\end{equation}
\end{lem}

\begin{proof}
By \eqref{eq:kernel-components}, injectivity is equivalent to $\m^nA_v\cap\m^{n+1}B_v=\m^{n+1}A_v$
for all $v,n$. An induction on $n$ gives the equivalent condition
\eqref{eq:strictness}.
\end{proof}

\begin{cor}\label{cor:strict}
Under the ring and module assumptions of
Theorem~\ref{thm:small-kernel}, suppose that $\cG_E$ is
equidimensional and $A_v\cap\m^nB_v=\m^nA_v$ for all $v,n\geq0.$
Then the following are equivalent:
\begin{enumerate}
\item[(i)] $U$ is a reduction of $E$;
\item[(ii)] $\br_i(U)=\br_i(E)$ for all $0\leq i\leq r$;
\item[(iii)] $\sum_{i=0}^r\br_i(U)=\sum_{i=0}^r\br_i(E)$.
\end{enumerate}
\end{cor}

\begin{proof}
By Lemma~\ref{lem:strictness}, the strictness condition gives $K=0$.
The conclusion follows from Theorem~\ref{thm:small-kernel}.
\end{proof}

\begin{cor}\label{cor:equigenerated}
Let $S$ be a finitely generated standard graded domain over a field, with $\dim S=d\geq1$, and set $\mathfrak n=S_+$ and $R=S_{\mathfrak n}$. Let $N\subseteq M\subseteq S^p$ be nonzero homogeneous modules of the same rank, both generated in the same degree $\delta\geq0$, with the basis of $S^p$ in degree zero. Set $U=N_{\mathfrak n}$ and $E=M_{\mathfrak n}$. Then $U$ is a reduction of $E$ if and only if $\br_i(U)=\br_i(E)$ for all $i$. 
\end{cor}

\begin{proof}
Inside $\Sym_S^v(S^p)$, the modules $M^v$ and $N^v$ are generated in degree $v\delta$. Since $S$ is standard graded, $\mathfrak n^nM^v=(M^v)_{\geq v\delta+n}$ and $\mathfrak n^nN^v=(N^v)_{\geq v\delta+n}$. Hence $N^v\cap\mathfrak n^nM^v=\mathfrak n^nN^v$, and after localization at $\mathfrak n$ we obtain the strictness condition \eqref{eq:strictness} for $A=\cR(U)$ and $B=\cR(E)$. Moreover, $\mathfrak n^nM^v/\mathfrak n^{n+1}M^v\simeq(M^v)_{v\delta+n}$, and these vector spaces are unchanged by localization at $\mathfrak n$. Thus $\cG_E$ identifies with $\cR_S(M)$ up to the regrading in which bidegree $(v,n)$ corresponds to degree $v\delta+n$ in the $v$-th Rees component. Since $\cR_S(M)$ embeds in $S[T_1,\ldots,T_p]$, it is a domain and hence equidimensional. Finally, $R$ is an excellent local domain, so the conclusion follows from Corollary~\ref{cor:strict}.
\end{proof}

\section{Relative multiplicities and reduction}
\label{sec:relative}

We now turn to relative polar multiplicities, which retain information
about the inclusion $U\subseteq E$ that is not visible in the two
absolute polar sequences.

\begin{nota}\label{not:relative}
For $A=\cR(U)\subseteq B=\cR(E)$, set $G=\gr_{A_+B}(B)$, $J=B_1G$, and $\br_i(U,E)=m_r^i(G).$
Here $G$ is endowed with the internal grading inherited from $B$, as in
\cite[Section~4]{CidRuiz2024}. The filtration degree is not added to
the internal degree. The notation $B_1G$ refers to the ideal generated
by the images in $G$ of the degree-one elements of $B$. There is a
graded exact sequence
\begin{equation}\label{eq:relative-exact}
0\longrightarrow J\longrightarrow G\longrightarrow A
\longrightarrow0.
\end{equation}
\end{nota}

The relative polar sequence provides a direct reduction criterion: its
agreement with the absolute polar sequence of $U$ is already equivalent
to reduction.

\begin{prop}
\label{prop:relative-absolute}
Let $R$ be an equidimensional universally catenary Noetherian local
ring of dimension $d\geq1$, and let $U\subseteq E$ be finite modules
of the same rank $e\geq1$. Let the Rees algebras be defined by
quotienting the symmetric algebras by torsion with respect to the
nonzerodivisors of $R$, and set $r=d+e-1$. Then the following
conditions are equivalent:
\begin{enumerate}
\item[(i)] $U$ is a reduction of $E$;
\item[(ii)] $\br_i(U,E)=\br_i(U)$ for every $0\leq i\leq r$;
\item[(iii)] $\sum_{i=0}^r(\br_i(U,E)-\br_i(U))=0$.
\end{enumerate}
Moreover, $\br_i(U,E)-\br_i(U)\geq0$ for every $0\leq i\leq r$.
\end{prop}

\begin{proof}
We first verify the equidimensionality required to apply
\cite[Theorem~4.4(iv)]{CidRuiz2024}. Choose a nonzerodivisor
$c\in R$ such that $E_c$ is free of rank $e$. Then
$B=\cR(E)$ embeds in $B_c\simeq R_c[T_1,\ldots,T_e]$. Hence every
minimal prime $\mathfrak q$ of $B$ contracts to a minimal prime
$\mathfrak p$ of $R$, and the dimension formula gives $\dim(B/\mathfrak q)
=
\dim(R/\mathfrak p)+e
=
d+e.$
Thus $B$ is equidimensional. Since $B$ is of finite type over the
universally catenary ring $R$, it is also universally catenary.

Because $U$ and $E$ have the same rank, the inclusion
$U\subseteq E$ becomes an isomorphism after tensoring with the total
quotient ring of $R$. Thus, whenever $A=\cR(U)\subseteq B=\cR(E)$ is
finite, the extension is birational. Therefore
\cite[Theorem~4.4(iv)]{CidRuiz2024} gives that (i) holds if and only
if (ii) holds. Applying additivity of polar multiplicities to
\eqref{eq:relative-exact} yields
\begin{equation}\label{eq:relative-defect}
\br_i(U,E)-\br_i(U)=m_r^i(J)\geq0
\end{equation}
for every $i$. Hence all these differences vanish if and only if their
sum vanishes, proving that (ii) holds if and only if (iii) holds.
\end{proof}

We next recall the classical Buchsbaum--Rim multiplicity for a finite-colength submodule of a free module, introduced in \cite{BuchsbaumRim1964} (see also \cite[Section~16.5]{HunekeSwanson}). Unlike the intrinsic polar multiplicities considered above, this invariant is defined relative to the ambient free module.

\begin{defn}
\label{def:classical-br}
Let $(R,\m)$ be a Noetherian local ring of positive dimension $d$,
let $F$ be free of rank $e\geq1$, and let $V\subseteq F$ have
finite colength. In this definition, $V^v$ means its image in
$\Sym_R^v(F)$. The classical multiplicity $e_{\mathrm{BR}}(V,F)$
is specified by
\[
\lambda_R(\Sym_R^v(F)/V^v)
=\frac{e_{\mathrm{BR}}(V,F)}{(d+e-1)!}v^{d+e-1}
 +O(v^{d+e-2}).
\]
We include $V=F$, with multiplicity zero.
\end{defn}

\begin{prop}
\label{prop:classical-reduction}
Let $R$ be an equidimensional universally catenary Noetherian local
ring of positive dimension, and let $U\subseteq E\subseteq F$, where
$F$ is finite free and $\lambda_R(F/U)<\infty$. Then $U$ is a
reduction of $E$ if and only if
$e_{\mathrm{BR}}(U,F)=e_{\mathrm{BR}}(E,F)$.
\end{prop}

\begin{proof}
This is \cite[Corollary~5.5]{SUV2001}, with the powers taken inside
$\Sym_R(F)$. Equivalently, equality of the two Buchsbaum--Rim
multiplicities is the vanishing of the leading term of
$\lambda_R(E^v/U^v)$, since
\[
\lambda_R(E^v/U^v)
=
\lambda_R(\Sym_R^v(F)/U^v)
-
\lambda_R(\Sym_R^v(F)/E^v).
\]
The cited reduction theorem characterizes precisely this degree drop.
\end{proof}

\begin{nota}\label{not:pair-coefficient}
Let $U\subseteq E\subseteq F=R^e$ have rank $e$, where $R$ is an
excellent local domain of dimension $d\geq1$, and assume
$\lambda_R(E/U)<\infty$. Set $r=d+e-1$. We denote by
$e_{\mathrm{pair}}(U,E)$ the relative Buchsbaum--Rim multiplicity
defined by
\begin{equation}\label{eq:pair-coefficient}
\lambda_R(E^v/U^v)
=
\frac{e_{\mathrm{pair}}(U,E)}{r!}v^r+O(v^{r-1}).
\end{equation}
The eventual polynomiality follows from
\cite[Theorem~16.5.6]{HunekeSwanson}. In the complex analytic
finite-colength setting, this is the multiplicity $e(U,E)$ used by
Gaffney (see \cite{KT1994}). We use the notation
$e_{\mathrm{pair}}(U,E)$ to distinguish it from the absolute polar
multiplicities.
\end{nota}

\begin{prop}
\label{prop:finite-relative-length}
Under Notation~\ref{not:pair-coefficient},
\begin{equation}\label{eq:finite-relative-polar}
\br_0(U,E)=e_{\mathrm{pair}}(U,E),\qquad
\br_i(U,E)=\br_i(U)\quad(1\leq i\leq r).
\end{equation}
In particular, $U$ is a reduction of $E$ if and only if
$\br_0(U,E)=0$. 
\end{prop}

\begin{proof}
Let $A=\cR(U)$, $B=\cR(E)$, $G=\gr_{A_+B}(B)$, and $J=B_1G$.
Since $\lambda_R(E/U)<\infty$, we have $A_{\mathfrak p}=B_{\mathfrak p}$
for every $\mathfrak p\neq\m$, and hence $J_{\mathfrak p}=0$.
As $J$ is a finite $G$-module, it follows that $\m^qJ=0$ for some
$q\geq1$. For each $v$, the internal grading gives
\[
J_v\simeq\bigoplus_{j=1}^{v}
\frac{A_{j-1}B_{v-j+1}}{A_jB_{v-j}},
\]
and therefore one has $\lambda_R(J_v)
=
\lambda_R(B_v/A_v)
=
\lambda_R(E^v/U^v).$
Since $\m^qJ=0$, for every $n\geq q-1$ and every $v$, $\lambda_R(J_v/\m^{n+1}J_v)=\lambda_R(J_v).$
Thus \eqref{eq:pair-coefficient} gives
$m_r^0(J)=e_{\mathrm{pair}}(U,E)$ and $m_r^i(J)=0$ for $i>0$. Applying polar additivity to
$0\to J\to G\to A\to0$ and using $\br_0(U)=0$, we obtain
\eqref{eq:finite-relative-polar}. The final assertion follows from
Proposition~\ref{prop:relative-absolute}.
\end{proof}

For a chain $U\subseteq E\subseteq P$ of modules of the same rank with finite relative lengths, the exact sequence of degree-$v$ quotients gives
\begin{equation}\label{eq:pair-additivity}
e_{\mathrm{pair}}(U,P)
=
e_{\mathrm{pair}}(U,E)
+
e_{\mathrm{pair}}(E,P).
\end{equation}
This is the finite-relative-length form of the additivity used in Gaffney's work.

\begin{prop}
\label{prop:relative-classical}
Let $R$ be an excellent local domain of dimension $d\geq1$, and let
$U\subseteq E\subseteq F=R^e$ have finite colength. Then
\[
\br_0(U,E)=e_{\mathrm{BR}}(U,F)-e_{\mathrm{BR}}(E,F),
\qquad
\br_i(U,E)=\br_i(U)\quad(1\leq i\leq r).
\]
In particular, $\br_0(U,F)=e_{\mathrm{BR}}(U,F)$, and $U$ is a
reduction of $E$ if and only if $\br_0(U,E)=0$.
\end{prop}

\begin{proof}
In the finite-colength case, for every $v\geq1$,
\begin{equation}\label{eq:relative-length}
\lambda_R(E^v/U^v)
=
\lambda_R(\Sym_R^v(F)/U^v)
-
\lambda_R(\Sym_R^v(F)/E^v).
\end{equation}
Comparing leading coefficients gives $e_{\mathrm{pair}}(U,E)
=
e_{\mathrm{BR}}(U,F)-e_{\mathrm{BR}}(E,F).$
The result now follows from
Proposition~\ref{prop:finite-relative-length}.
\end{proof}

\begin{rem}\label{rem:false-identification}
The identity $\br_0(U,F)=e_{\mathrm{BR}}(U,F)$ is relative to the
ambient free module and should not be confused with the absolute
invariant $\br_0(U)=0$. Thus the classical Buchsbaum--Rim criterion
does not yield an absolute polar criterion.
\end{rem}

For comparison with \cite[Section~16.5]{HunekeSwanson}, their
Buchsbaum--Rim invariants correspond, under the finite-colength
hypotheses, to
\[
\br^{\mathrm{HS}}(U)=e_{\mathrm{BR}}(U,F)=\br_0(U,F),
\qquad
\br^{\mathrm{HS}}(U,E)
=
e_{\mathrm{BR}}(U,F)-e_{\mathrm{BR}}(E,F)
=
\br_0(U,E).
\]
In particular, their reduction criterion for finite-colength pairs is
a relative statement: vanishing of $\br^{\mathrm{HS}}(U,E)$ detects
reduction, while the absolute invariant $\br_0(U)$ remains zero.

\section{Ideal multiples and relative multiplicities}
\label{sec:ideal-multiples}

We first compare reductions and polar multiplicities of ideals with
those of their free-module multiples. We then study the effect of
maximal-ideal twists, which will be used in the counterexamples of
Section~\ref{sec:counterexamples}.

Return to the standing assumptions: $R$ is an excellent local domain
of dimension $d\geq1$, and let $F=R^e$, with $e\geq1$.

\begin{prop}
\label{prop:ideal-transfer}
Let $0\neq I\subseteq J$ be ideals of $R$. Then, for all $v,n\geq0$,
\begin{equation}\label{eq:ideal-transfer-hilbert}
H_{IF}(v,n)
=
\binom{v+e-1}{e-1}H_I(v,n).
\end{equation}
Consequently,
\begin{equation}\label{eq:ideal-transfer-polar}
\br_i(IF)
=
\binom{d+e-1-i}{e-1}\br_i(I)
\quad(0\leq i\leq d),
\end{equation}
and $\br_i(IF)=0$ for $d<i\leq d+e-1$. Moreover, $IF$ is a reduction
of $JF$ if and only if $I$ is a reduction of $J$, and the polar
sequences of $IF$ and $JF$ are equal if and only if the polar
sequences of $I$ and $J$ are equal.
\end{prop}

\begin{proof}
Using the monomial basis of $\Sym_R^v(F)$, we have
\[
(IF)^v
=
I^v\Sym_R^v(F)
\simeq
(I^v)^{\oplus\binom{v+e-1}{e-1}}.
\]
This identification is compatible with multiplication by powers of
$\m$, and hence gives \eqref{eq:ideal-transfer-hilbert}. Since the
leading term of $\binom{v+e-1}{e-1}$ is $v^{e-1}/(e-1)!$, comparison
with \eqref{eq:polar-polynomial} yields
\[
\br_i(IF)
=
\binom{d+e-1-i}{e-1}\br_i(I)
\quad(0\leq i\leq d),
\]
while the remaining polar multiplicities vanish. Since all the
binomial factors are positive, equality of the polar sequences of
$IF$ and $JF$ is equivalent to equality of those of $I$ and $J$.
Hence, for every $v\geq0$, $(JF)^{v+1}
=
J^{v+1}\Sym_R^{v+1}(F),$ and $(IF)(JF)^v
=
IJ^v\Sym_R^{v+1}(F).$
Thus $(JF)^{v+1}=(IF)(JF)^v$ if and only if
$J^{v+1}=IJ^v$. Hence $IF$ is a reduction of $JF$ if and only if
$I$ is a reduction of $J$.
\end{proof}

\begin{prop}
\label{prop:ideal-transfer-classical}
Let $I\subseteq J$ be $\m$-primary ideals, and let $e(I)$ denote the Hilbert--Samuel multiplicity of $I$. Then
\[
e_{\mathrm{BR}}(IF,F)=\binom{d+e-1}{e-1}e(I),
\qquad
\br_0(IF,JF)=\binom{d+e-1}{e-1}\bigl(e(I)-e(J)\bigr).
\]
Moreover, $\br_i(IF,JF)=\br_i(IF)$ for every $i\geq1$.
\end{prop}

\begin{proof}
For every $v\geq1$,
\[
\lambda_R\bigl(\Sym_R^v(F)/(IF)^v\bigr)
=
\binom{v+e-1}{e-1}\lambda_R(R/I^v).
\]
Comparing leading terms gives
\[
e_{\mathrm{BR}}(IF,F)
=
\binom{d+e-1}{e-1}e(I).
\]
The remaining assertions follow from
Proposition~\ref{prop:relative-classical}.
\end{proof}

For every $a\geq0$, the inclusion $\m^aM\subseteq M$ has finite relative length, even when $F/M$ does not. The next result shows that the absolute polar sequence is encoded by the corresponding relative Buchsbaum--Rim multiplicities.

\begin{thm}
\label{thm:twist-polynomial}
Let $R$ and $M$ be as above, and set $r=d+e-1$. Then, for every integer $a\geq0$,
\begin{equation}\label{eq:twist-polynomial}
e_{\mathrm{pair}}(\m^aM,M)
=
\sum_{i=1}^{d}\binom{r}{i}a^i\br_i(M).
\end{equation}
If $\lambda_R(F/M)<\infty$, then
\begin{equation}\label{eq:twist-classical}
e_{\mathrm{BR}}(\m^aM,F)
=
e_{\mathrm{BR}}(M,F)
+
\sum_{i=1}^{d}\binom{r}{i}a^i\br_i(M).
\end{equation}
In particular, the values $e_{\mathrm{pair}}(\m^aM,M)$ for $a=1,\ldots,d$ determine the polar sequence of $M$ once the rank and dimension are fixed. If $R$ is fixed, then the values for $a=1,\ldots,d-1$ already suffice, since $\br_d(M)=e_{\m}(R)$. 
\end{thm}

\begin{proof}
Assume first that $a\geq1$. Since $(\m^aM)^v=\m^{av}M^v$, for every $v\geq1$ we have
\begin{equation}\label{eq:twist-diagonal}
\lambda_R\bigl(M^v/(\m^aM)^v\bigr)
=
H_M(v,av-1).
\end{equation}
For $v\gg0$, the right-hand side lies in the polynomial range of $H_M$. Substituting $n=av-1$ into \eqref{eq:polar-polynomial}, the coefficient of $v^r$ is
\[
\sum_{i=0}^r\frac{a^i\br_i(M)}{(r-i)!i!}.
\]
Hence, by the definition of $e_{\mathrm{pair}}(\m^aM,M)$,
\[
e_{\mathrm{pair}}(\m^aM,M)
=
\sum_{i=0}^r\binom{r}{i}a^i\br_i(M).
\]
Proposition~\ref{prop:fixed-entries} gives $\br_0(M)=0$ and $\br_i(M)=0$ for $i>d$, yielding \eqref{eq:twist-polynomial}. The case $a=0$ is immediate. In addition, if $\lambda_R(F/M)<\infty$, additivity of relative Buchsbaum--Rim multiplicity gives
\[
e_{\mathrm{BR}}(\m^aM,F)
=
e_{\mathrm{BR}}(M,F)
+
e_{\mathrm{pair}}(\m^aM,M),
\]
and \eqref{eq:twist-classical} follows. The right-hand side of \eqref{eq:twist-polynomial} is a polynomial in $a$ of degree at most $d$ with zero constant term. Its values at $d$ distinct nonzero integers determine its coefficients and hence $\br_1(M),\ldots,\br_d(M)$. If $R$ is fixed, the coefficient involving $\br_d(M)=e_{\m}(R)$ is already known, so $d-1$ values suffice.
\end{proof}

\section{Counterexamples}
\label{sec:counterexamples}

We now show that the failure of the converse in
Question~\ref{question:main} is not accidental. On regular surfaces,
the determinant formula produces broad families of distinct modules
with identical intrinsic polar data. 

\begin{prop}
\label{prop:closed-counterexamples}
Let $(R,\m,k)$ be an excellent regular local ring of dimension two
with $k$ infinite. Let $U\subsetneq E\subseteq F=R^e$ have finite
colength, and assume that $U$ is integrally closed in $F$. If
$o(U)=o(E)$, then $U$ and $E$ have the same polar sequence, but
$U$ is not a reduction of $E$.
\end{prop}

\begin{proof}
By Theorem~\ref{thm:surface-polars}, the equality $o(U)=o(E)$ implies
that $U$ and $E$ have the same polar sequence. If $U$ were a reduction
of $E$, then $E\subseteq\overline U=U$, contradicting
$U\subsetneq E$. 
\end{proof}

A particularly  source of counterexamples comes from complete ideals of the same order.

\begin{cor}
\label{cor:complete-ideal-counterexamples}
Let $I\subsetneq J$ be $\m$-primary ideals with $I$ integrally closed and
$\operatorname{ord}_{\m}I=\operatorname{ord}_{\m}J$. Then
$I\subsetneq J$ is a counterexample to Question~\ref{question:main}.
More generally, for every $e\geq1$ and $F=R^e$, the inclusion
$IF\subsetneq JF$ is also a counterexample, with common polar sequence
\[
(0,e\operatorname{ord}_{\m}I,1,0,\ldots,0).
\]
\end{cor}

\begin{proof}
The rank-one case follows from
Proposition~\ref{prop:closed-counterexamples}. The general case follows
from Proposition~\ref{prop:ideal-transfer}.
\end{proof}

The next family strengthens the preceding counterexamples: the entire intrinsic Hilbert function remains unchanged, although reduction fails.

\begin{prop}\label{prop:plane-example}
Let $R=k[[x,y]]$, $\m=(x,y)$, and $I_q=(x,y^q)$ for $q\geq1$. Then, for every $v,n\geq0$,
\begin{equation}\label{eq:exact-hilbert}
H_{I_q}(v,n)=\binom{n+2}{2}+v(n+1).
\end{equation}
Consequently, $(\br_0(I_q),\br_1(I_q),\br_2(I_q))=(0,1,1).$
For every $q\geq2$, the ideal $I_q$ is not a reduction of $\m=I_1$.
\end{prop}

\begin{proof}
Since $x,y^q$ is a regular sequence, $I_q$ is of linear type, and hence $\cR(I_q)\simeq R[S,T]/(y^qS-xT).$
With respect to the $\m$-adic filtration on the coefficients,
$\gr_{\m}R\simeq k[X,Y]$ is a domain. Thus the initial ideal of
$(y^qS-xT)$ is generated by the initial form of its defining relation,
and therefore
\begin{equation}\label{eq:initial-rings}
\cG_{I_q}\simeq
\begin{cases}
k[X,Y,S,T]/(YS-XT),&q=1,\\
k[X,Y,S,T]/(XT),&q\geq2.
\end{cases}
\end{equation}
Give $S,T$ bidegree $(1,0)$ and $X,Y$ bidegree $(0,1)$. In either
case the defining equation is a nonzerodivisor of bidegree $(1,1)$,
so
\[
\Hilb_{\cG_{I_q}}(s,t)
=
\frac{1-st}{(1-s)^2(1-t)^2}.
\]
Since $H_{I_q}(v,n)$ is obtained by summing the graded pieces in the
second degree,
\[
\sum_{v,n\geq0}H_{I_q}(v,n)s^vt^n
=
\frac{1-st}{(1-s)^2(1-t)^3}.
\]
Taking the coefficient of $s^vt^n$ gives
\[
H_{I_q}(v,n)
=
(v+1)\binom{n+2}{2}
-
v\binom{n+1}{2}
=
\binom{n+2}{2}+v(n+1),
\]
proving \eqref{eq:exact-hilbert}. Its homogeneous part of total degree
two is $vn+n^2/2$, and comparison with
\eqref{eq:polar-polynomial} gives the stated polar sequence.
Hence, if $I_q$ were a reduction of $\m$ for some $q\geq2$, then
$\m^{v+1}=I_q\m^v$ for some $v$. Reducing modulo $(x)$ would give
$(y^{v+1})=(y^{v+q})$ in $k[[y]]$, a contradiction.
\end{proof}

\begin{rem}\label{rem:strength-example}
The ideals $I_q$ and $\m$ are $\m$-primary, have the same order, and
coincide away from $\m$, with $\lambda_R(\m/I_q)=q-1$. Thus even under
these hypotheses, equality of the entire intrinsic Hilbert function
does not imply reduction.
\end{rem}

Note that Proposition~\ref{prop:plane-example} shows that even equality of the full intrinsic Hilbert function does not detect reduction in rank one. The same phenomenon persists in every positive rank, while keeping finite colength and an explicitly computable polar sequence.

\begin{cor}\label{cor:all-ranks}
For $e\geq1$ and $q\geq2$, set $U=I_q\oplus R^{e-1}\subseteq E=\m\oplus R^{e-1}\subseteq R^e.$
Then $U$ is not a reduction of $E$, although
\[
\br_i(U)=\br_i(E)=
\begin{cases}
1,&i=1,2,\\
0,&0\leq i\leq e+1,\ i\neq1,2.
\end{cases}
\]
Moreover, both modules have finite colength in $R^e$.
\end{cor}

\begin{proof}
For $I=I_q$ or $\m$, one has
$\cR(I\oplus R^{e-1})\simeq\cR(I)[Z_1,\ldots,Z_{e-1}]$, where the
variables $Z_j$ have bidegree $(1,0)$. Using
Proposition~\ref{prop:plane-example}, the generating series of the
corresponding Hilbert function is
\[
\frac{1-st}{(1-s)^{e+1}(1-t)^3}.
\]
Hence
\[
H(v,n)
=
\binom{v+e-1}{e-1}\binom{n+2}{2}
+
\binom{v+e-1}{e}(n+1).
\]
Its homogeneous part of total degree $e+1$ is
\[
\frac{v^{e-1}n^2}{(e-1)!2!}
+
\frac{v^en}{e!},
\]
which gives the stated polar sequence. If $U$ were a reduction of $E$, setting $Z_1=\cdots=Z_{e-1}=0$
in a reduction identity would imply that $I_q$ is a reduction of
$\m$, contradicting Proposition~\ref{prop:plane-example}. Finally,
$\lambda_R(R^e/U)=q$ and $\lambda_R(R^e/E)=1$.
\end{proof}

It should be noted that Proposition~\ref{prop:plane-example} also shows that the strict dimension bound in Theorem~\ref{thm:small-kernel} is sharp: for this counterexample, the kernel of the associated graded map has dimension exactly $d+e$.

\begin{prop}\label{prop:explicit-kernel}
Let $A=\cR(I_q)\subseteq B=\cR(\m)$ with $q\geq2$, as in
Proposition~\ref{prop:plane-example}. Then the associated graded map is
\[
\frac{k[X,Y,S,T]}{(XT)}
\longrightarrow
\frac{k[X,Y,P,Q]}{(YP-XQ)},
\qquad
S\longmapsto P,\quad T\longmapsto0,
\]
and fixes $X,Y$. Its kernel is $K=(T)$ and its image is
$k[X,Y,P]$. In particular, $\dim K=3=d+e,$
so the strict inequality in Theorem~\ref{thm:small-kernel} cannot be
weakened to $\dim K\leq d+e$.
\end{prop}

\begin{proof}
The generators $S$ and $T$ correspond respectively to $xt$ and
$y^qt$. Hence $S$ maps to $P$, whereas
$y^qt=y^{q-1}(yt)\in\m B_1$, so $T$ maps to zero in the associated
graded algebra. Thus $(T)$ is contained in the kernel. On the other hand,
$X,Y,P$ are algebraically independent in
$k[X,Y,P,Q]/(YP-XQ)$, and therefore the induced map
$k[X,Y,S]\to k[X,Y,P]$ is an isomorphism. Hence $K=(T)$ and the image
is $k[X,Y,P]$. Therefore, in $k[X,Y,S,T]/(XT)$ one has
$\operatorname{Ann}(T)=(X)$. So $K\simeq k[Y,S,T](-1,0),$
and hence $\dim K=3$. Since here $d=2$ and $e=1$, this is exactly
$d+e$.
\end{proof}

\begin{lem}\label{lem:staircase-count}
Let $R=k[[x,y]]$, and let $L=(x^a,x^{a-1}y^{b_1},\ldots,y^{b_a})$, where $0=b_0<b_1<\cdots<b_a,$
where $a\geq0$, with $L=R$ when $a=0$. Then, for every $n\geq0$,
\begin{equation}\label{eq:staircase-count}
\lambda_R(L/\m^{n+1}L)=\binom{n+2}{2}+a(n+1).
\end{equation}
\end{lem}

\begin{proof}
Set $t=n+1$. For $0\leq h<a$, the least $y$-exponent of a monomial
in $L$ with $x$-exponent $h$ is $b_{a-h}$. In $\m^tL$, the
corresponding least exponent is $b_{a-h}+t$. Indeed, passing to a
generator with smaller $x$-exponent increases its $y$-exponent by at
least the number of steps, since $b_{i+j}\geq b_i+j$. Thus these $a$
columns contribute $at$ to the difference of colengths.

For $a\leq h<a+t$, the least $y$-exponents in $L$ and $\m^tL$ are,
respectively, $0$ and $a+t-h$, while for $h\geq a+t$ both are zero.
Hence the remaining contribution is
\[
t+(t-1)+\cdots+1=\binom{t+1}{2}.
\]
Therefore
\[
\lambda_R(L/\m^tL)
=
at+\binom{t+1}{2}
=
\binom{n+2}{2}+a(n+1),
\]
as claimed.
\end{proof}

The next family produces counterexamples to Question~\ref{question:main} with substantially stronger properties. The intrinsic Hilbert function, and hence the entire polar sequence, is independent of the parameters $p$ and $q$, while the ideals themselves may vary through distinct normal $\m$-primary ideals. 

\begin{thm}
\label{thm:two-face-family}
In $R=k[[x,y]]$, let $K_{p,q}=(x^2,xy^p,y^{p+q})$,  where $1\leq p\leq q.$
Then every power of $K_{p,q}$ is integrally closed and, for all
$v,n\geq0$,
\begin{align}
H_{K_{p,q}}(v,n)
 &=\binom{n+2}{2}+2v(n+1),
 \label{eq:two-face-hilbert}\\
\lambda_R(R/K_{p,q}^v)
 &=\frac{3p+q}{2}v^2+\frac{p+q}{2}v.
 \label{eq:two-face-colength}
\end{align}
Moreover, for $F=R^e$,
\begin{align}
H_{K_{p,q}F}(v,n)
 &=\binom{v+e-1}{e-1}
 \left[\binom{n+2}{2}+2v(n+1)\right],
 \label{eq:two-face-module-hilbert}\\
\bigl(\br_i(K_{p,q}F)\bigr)_{i=0}^{e+1}
 &=(0,2e,1,0,\ldots,0),
 \label{eq:two-face-module-polar}\\
e_{\mathrm{BR}}(K_{p,q}F,F)
 &=\binom{e+1}{2}(3p+q).
 \label{eq:two-face-module-classical}
\end{align}
The module $K_{p,q}F$ is integrally closed in $F$ and has no nonzero
free direct summand.
\end{thm}

\begin{proof}
The minimal monomial generators of $K_{p,q}^v$ are $x^{2v-j}y^{c_j}$ for $0\leq j\leq2v,$ where
\[
c_j=
\begin{cases}
pj,&0\leq j\leq v,\\
qj+(p-q)v,&v\leq j\leq2v.
\end{cases}
\]
Indeed, for a fixed loss $j$ in the $x$-exponent, the least possible
$y$-exponent is obtained by first replacing copies of $x^2$ by
$xy^p$ and then replacing copies of $xy^p$ by $y^{p+q}$; the
corresponding costs are $p$ and $q\geq p$. Since the successive
differences of the sequence $(c_j)$ are $p$ and $q$,
Lemma~\ref{lem:staircase-count} gives
\eqref{eq:two-face-hilbert}. Moreover,
\[
\lambda_R(R/K_{p,q}^v)
=\sum_{j=1}^{2v}c_j
=\frac{3p+q}{2}v^2+\frac{p+q}{2}v,
\]
proving \eqref{eq:two-face-colength}. To prove integral closedness, define valuations on the fraction field
of $R$ by $\nu_p(x)=p$, $\nu_p(y)=1$,
$\nu_q(x)=q$, $\nu_q(y)=1$.
For every nonzero power series, each valuation is the minimum of the
weights of its monomials. The above description of the generators
shows that
\begin{equation}\label{eq:two-face-valuations}
K_{p,q}^v
=
\{f\in R:\nu_p(f)\geq2pv,\ 
           \nu_q(f)\geq(p+q)v\}.
\end{equation}
Indeed, for a monomial with $x$-exponent $h=2v-j$, the least allowed
$y$-exponent is
\[
\max\{0,\,2pv-ph,\,(p+q)v-qh\}=c_j.
\]
Each valuation inequality is preserved under integral dependence over
an ideal. Hence the right-hand side of
\eqref{eq:two-face-valuations} is integrally closed, and so every
$K_{p,q}^v$ is integrally closed. When $q>p$, the two inequalities
correspond to two distinct compact faces of the Newton boundary.
Now, Equations \eqref{eq:two-face-module-hilbert}--
\eqref{eq:two-face-module-classical} follow from
Propositions~\ref{prop:ideal-transfer} and
\ref{prop:ideal-transfer-classical}. To see that $K_{p,q}F$ is
integrally closed, let $z=(z_1,\ldots,z_e)\in F$ be integral over
$K_{p,q}F$. Applying each coordinate projection to the corresponding
integral equation in the symmetric algebra shows that every $z_i$ is
integral over $K_{p,q}$. Since $K_{p,q}$ is integrally closed,
$z_i\in K_{p,q}$ for every $i$, and hence $z\in K_{p,q}F$.

Let $L$ be any $\m$-primary ideal of $R$. Since $R$ is a
two-dimensional regular local ring, every homomorphism $L\to R$ is
multiplication by an element of $R$. Indeed, after tensoring with the
fraction field it is multiplication by some scalar $c$, and
$L_{\mathfrak p}=R_{\mathfrak p}$ for every height-one prime
$\mathfrak p$, so $c\in R_{\mathfrak p}$ for all such $\mathfrak p$.
As $R$ is normal,
$c\in\bigcap_{\operatorname{ht}\mathfrak p=1}R_{\mathfrak p}=R$.
Thus every map $L^{\oplus e}\to R$ has image contained in
$L\subseteq\m$ and cannot be surjective. Taking $L=K_{p,q}$ shows
that $K_{p,q}F$ has no nonzero free direct summand.
\end{proof}

 The following rank-two example makes the failure completely concrete, even for integrally closed modules with identical intrinsic Hilbert functions and no nonzero free direct summands.

\begin{ex}\label{ex:two-face-rank-two}
Let $R=k[[x,y]]$, $F=R^2$, and set $I=(x^2,xy^2,y^5)=K_{2,3}$, $J=(x^2,xy,y^3)=K_{1,2}$, and
$U=IF\subsetneq E=JF.$
Both $I$ and $J$ are normal ideals with two distinct compact Newton
faces. Consequently, $U$ and $E$ are integrally closed, have six
minimal generators, and have no nonzero free direct summand. Nevertheless,
\[
H_U(v,n)=H_E(v,n)
=(v+1)\left[\binom{n+2}{2}+2v(n+1)\right]
\]
for every $v,n\geq0$, and hence both modules have polar sequence
$(0,4,1,0)$. On the other hand,
\[
\begin{array}{c|cc}
 &U&E\\ \hline
\lambda_R(F/(-))&14&8\\
e_{\mathrm{BR}}(-,F)&27&15
\end{array}
\qquad
\bigl(\br_i(U,E)\bigr)_{i=0}^{3}=(12,4,1,0).
\]
Also, $I$ is not a reduction of $J$: if
$J^{v+1}=IJ^v$ for some $v$, then modulo $(x)$ one would obtain
$(y^{3v+3})=(y^{3v+5})$ in $k[[y]]$, a contradiction.
So Proposition~\ref{prop:ideal-transfer} therefore shows that $U$ is not
a reduction of $E$.
\end{ex}

The preceding family gives not only counterexamples to
Question~\ref{question:main}, but counterexamples with arbitrarily
large relative defect while all intrinsic Hilbert data remain fixed.

\begin{cor}
\label{cor:unbounded-two-face}
Fix $e\geq1$ and $p\geq1$. For integers $q'>q>p$, set
$U=K_{p,q'}F$ and $E=K_{p,q}F$. Then
$U\subsetneq E$, $U$ is not a reduction of $E$, and
\[
H_U=H_E,\qquad
\br_0(U,E)=\binom{e+1}{2}(q'-q).
\]
In particular, $\br_0(U,E)$ is unbounded among integrally closed
modules of fixed rank, with the same intrinsic Hilbert function and
the same number $3e$ of minimal generators. Moreover, both endpoint
ideals are normal and have two distinct compact Newton faces.
\end{cor}

\begin{proof}
The inclusion $K_{p,q'}\subsetneq K_{p,q}$ is immediate from the
generators and is strict modulo $(x)$. Since $K_{p,q'}$ is integrally
closed, it cannot be a reduction of the strictly larger ideal
$K_{p,q}$. Proposition~\ref{prop:ideal-transfer} therefore shows that
$U$ is not a reduction of $E$.

By Theorem~\ref{thm:two-face-family}, the intrinsic Hilbert function
of $K_{p,q}F$ is independent of $q$, so $H_U=H_E$. Finally,
Proposition~\ref{prop:relative-classical} and
\eqref{eq:two-face-module-classical} give
\[
\br_0(U,E)
=
e_{\mathrm{BR}}(U,F)-e_{\mathrm{BR}}(E,F)
=
\binom{e+1}{2}(q'-q).
\]
\end{proof}

\begin{rem}
The results of this section give counterexamples to
Question~\ref{question:main} under several additional hypotheses.
More precisely, the converse fails even for finite-colength
submodules of a free module over a two-dimensional regular local
ring in each of the following situations:
\begin{enumerate}
\item the smaller module is integrally closed and the two modules
have the same determinant order (Proposition~\ref{prop:closed-counterexamples});

\item both endpoints are integrally closed, in which case their
entire intrinsic Hilbert functions may coincide (Proposition~\ref{prop:closed-counterexamples} and
Proposition~\ref{prop:plane-example});

\item the modules have arbitrary fixed positive rank and the same
absolute polar sequence (Corollary~\ref{cor:all-ranks});

\item both endpoints arise from normal $\m$-primary ideals whose
Newton boundaries have two distinct compact faces (Theorem~\ref{thm:two-face-family} and
Example~\ref{ex:two-face-rank-two});

\item both modules are integrally closed, have no nonzero free direct
summand, and may have the same rank, the same number of minimal
generators, and the same intrinsic Hilbert function, while the
relative defect $\br_0(U,E)$ is arbitrarily large (Theorem~\ref{thm:two-face-family} and
Corollary~\ref{cor:unbounded-two-face}).

\end{enumerate}

\end{rem}
\section{Acknowledgments}
This work was developed during the postdoctoral fellowship of the first author at the Instituto de Ci\^encias Matem\'aticas e de Computa\c{c}\~ao (ICMC), Universidade de S\~ao Paulo (USP) - 2026. He gratefully acknowledges financial support from FAPESP, grant 2025/20830-5. The first author also gratefully acknowledges the Universidade Tecnol\'ogica Federal do Paran\'a -- Campus Guarapuava for the opportunity to undertake his postdoctoral research. The second author thanks CNPq-Brazil, grant 304851/2025-6 and -FAPESP-Brazil, grant 2025/21618-0.
OpenAI’s ChatGPT was used during the preparation of this work as an auxiliary generative-AI tool for  discussion, checks of algebraic manipulations and examples, comparison with cited literature, and assistance with English exposition and LaTeX. All material used in the final manuscript was checked by the authors against the relevant arguments and sources. The authors assume full responsibility for the final text.

\end{document}